\documentclass[reqno]{article}
\usepackage{amsmath, amsfonts, amscd, eucal}
\usepackage{latexsym}
\usepackage{amssymb}
\newcommand{\be}{\begin{equation}}
      \newcommand{\ee}{\end{equation}}
      \newcommand{\ba}{\begin{eqnarray}}
       \newcommand{\ea}{\end{eqnarray}}
\newcommand{\ban}{\begin{eqnarray*}}
\newcommand{\ean}{\end{eqnarray*}}

\newcommand{\lp}{\langle}
\newcommand{\rp}{\rangle}
\newcommand{\ra}{\rightarrow}

 \renewcommand{\o}[2]{\frac{#1}{#2}}
\newcommand{\hf}{\o{1}{2}}

 \newcommand{\qed}{\hspace*{\fill}\rule{3mm}{3mm}\quad \vspace{.2cm}}
 \newcommand{\Pf}{\noindent {\bf Proof:} }
 \newcommand{\Rk}{\noindent {\bf Remark} }
\newcommand{\Ex}{\noindent {\bf Example:} }

\newcommand{\sect}[1]{\section{#1} \setcounter{equation}{0}}

\newtheorem{theo}{Theorem}[section]

\begin{document}
\newtheorem{defn}[theo]{Definition}
\newtheorem{ques}[theo]{Question}
\newtheorem{lem}[theo]{Lemma}
\newtheorem{prop}[theo]{Proposition}
\newtheorem{coro}[theo]{Corollary}
\newtheorem{ex}[theo]{Example}
\newtheorem{note}[theo]{Note}
\newtheorem{conj}[theo]{Conjecture}

\title{On the Stability of K\"ahler-Einstein Metrics}
\author{Xianzhe Dai\thanks {Math Dept, UCSB, Santa Barbara, CA 93106 \tt{Email:
dai@math.ucsb.edu}. Partially
supported by NSF Grant \# DMS-0405890 } \and Xiaodong Wang\thanks{Math Dept, Michigan State University,
East Lansing, MI 48824 \tt{Email:xwang@math.msu.edu}.}
\and Guofang Wei\thanks
{Math Dept. UCSB. \tt{Email: wei@math.ucsb.edu}. Partially
supported by NSF Grant \# DMS-0204187.}}
\maketitle

\begin{abstract}
Using spin$^c$ structure we prove that K\"ahler-Einstein metrics
with nonpositive scalar curvature are stable (in the direction of
changes in conformal structures) as the critical points of the
total scalar curvature functional. Moreover if all infinitesimal
complex deformation of the complex structure are integrable, then
the K\"ahler-Einstein metric is a local maximal of the Yamabe
invariant, and its volume is a local minimum among all metrics
with scalar curvature bigger or equal to the scalar curvature of
the K\"ahler-Einstein metric.
\end{abstract}

\newcommand{\bi}{\bar{i}}
\newcommand{\vol}{\mbox{vol}}
\newcommand{\diam}{\mbox{diam}}
\newcommand{\Ric}{\mbox{Ric}}
\newcommand{\Iso}{\mbox{Iso}}
\newcommand{\Hess}{\mbox{Hess}}
\newcommand{\divg}{\mbox{div}}
\newcommand\grd{\nabla}
\newcommand{\Dirac}{\mathcal{D}}
\newcommand{\M}{\mathcal{M}}
\newcommand{\lx}{L_X\raisebox{0.5ex}{$g$}}
\newcommand{\cuv}{\overset{\hspace{.5ex}\circ}{R}}
\newcommand{\Sr}{\mathcal{S}}
\newcommand{\ts}{\otimes}
\newcommand{\Mh}{\hat{M}}
\newcommand{\gh}{\hat{g}}

\def\operatorname#1{{\rm #1\,}}
\def\op{\operatorname}
\def\lam{\lambda}
\def\s{\sigma}
\def\ph{\phi}
\def\ps{\psi}
\def\e{\epsilon}
\def\d{\delta}
\def\D{\Delta}
\sect{Introduction}

Stability issue comes up naturally in variational problems. One of
the most important geometric variational problems is that of the
total scalar curvature functional. Following \cite[Page
132]{Besse} we call an Einstein metric stable if the second
variation of the total scalar curvature functional is non-positive
in the direction of changes in conformal structures (we have
weakened the notion by allowing kernels; see also \cite{ko2} where
stability is defined in terms of local maximum). By the well-known
formula, this is to say,
 \be \label{stability} \lp \grd^*\grd h-2\cuv h, h\rp \ge 0 \ee  for any trace-free and
divergence-free symmetric two tensor $h$. Here $\cuv h$ denotes
the natural action of the curvature tensor on the symmetric
tensors \cite{Besse}. The operator appearing in (\ref{stability})
is closely related to the Lichnerowicz Laplacian $\mathcal L_g$.
Indeed, one has \be \label{lichlap} \mathcal L_gh=\grd^*\grd
h-2\cuv h+  \Ric \circ h + h \circ \Ric. \ee The two thus coincide
for Ricci flat metrics.

In \cite{dww}, we studied the stability of compact Ricci flat
manifolds. An essential ingredient there is the use of spin structure and
parallel spinors. In fact, our result should really be viewed as
the stability result for compact Riemannian manifolds with nonzero
parallel spinor. By \cite{Wa}, \cite{H}, this class of manifolds essentially coincides with that of special holonomy, namely, the Calabi-Yau manifolds, hyperK\"ahler manifolds, spin(7) manifolds and $G_2$ manifolds.

In this paper, we use spin$^c$ structure to generalize our
previous result to manifolds with nonzero parallel spin$^c$
spinor. Since the existence of nonzero parallel spinor implies
that the metric is necessarily Ricci flat, our motivation here is
to extend our previous method to deal with nonzero scalar
curvature and we found spin$^c$ to be a good framework to work
with.

\begin{theo} If a compact Einstein manifold $(M, g)$ with nonpositive
scalar curvature admits a nonzero parallel spin$^c$ spinor,  then it is stable.
\end{theo}

As we mentioned, this generalizes the stability result in
\cite{dww}. Since a K\"ahler manifold with its canonical spin$^c$
structure has nonzero parallel spin$^c$ spinors, this implies

\begin{coro} A compact K\"ahler-Einstein manifold with
non-positive scalar curvature is stable.
\end{coro}
This also follows essentially from Koiso's work \cite{ko},
\cite{Besse}, although it does not seem to have been noticed
before. Our approach of using spin$^c$ structure is new and gives
more general result. A well known result in the same direction is
for compact Einstein manifolds with negative sectional curvature
\cite{ko1}, \cite{ye}, \cite{Besse}. In this case the manifold is
strictly stable in the sense that the operator $ \grd^*\grd -2\cuv
$ is in fact positive definite. In contrast, there are many
Einstein manifolds with positive scalar curvature which are
unstable \cite{ko2}, \cite{chi} (see also \cite{boehm}).

It turns out that manifolds admitting a nonzero parallel spin$^c$
spinor are more or less classified \cite{mo}. Namely a simply
connected manifold has a nonzero parallel spin$^c$ spinor if and
only if the manifold is the product of a K\"ahler manifold and a
manifold with parallel spinor.  Moreover, the spin$^c$ structure
is the product of the canonical spin$^c$ structure on the K\"ahler
manifold with the spin structure on the other factor.

For manifolds with nonzero parallel spin$^c$ spinor, we derive a
Bochner type formula relating the operator $ \grd^*\grd -2\cuv $
to the square of a twisted Dirac operator. The difference, which
is expressed in terms of the curvatures, can be shown to be
nonnegative under our assumption. In fact, we prove that the
operator $\grd^*\grd -2\cuv $ is positive semi-definite for
K\"ahler manifolds with nonpositive Ricci curvature. Our method
also proves that the Lichnerowicz Laplacian is positive
semi-definite for K\"ahler manifolds with nonnegative Ricci
curvature.

The operator $\grd^*\grd -2\cuv $ (or the Lichnerowicz Laplacian)
seems to have a knack for appearing in geometric variational
problems. Besides the total scalar curvature functional, or
equivalently, the Yamabe functional if one normalizes the volume,
\[ Y(g)=\frac{\int_M S_gdV_g}{{\rm Vol}(g)^{1-\frac{2}{n}}},\]
there is also the first eigenvalue $\lambda(g)$ of conformal
Laplacian considered in \cite{dww}, and the $L^{n/2}$ norm of
scalar curvature \cite{bcg}
\[ K(g)=\int_M|S_g|^{n/2}dV_g.  \]
Using these functionals, we can then deduce a number of
interesting consequences.

\begin{theo} 
Let $(N,g_0,J_0)$ be a compact K\"ahler-Einstein manifold with
nonpositive scalar curvature. Suppose all infinitesimal complex
deformations of $J_0$ are integrable. Then $g_0$ is a local
maximum of the Yamabe invariant.
\end{theo}

In the case of zero scalar curvature, the integrability  condition
is automatic by the Bogomolov-Tian-Todorov theorem \cite{Bo},
\cite{t}, \cite{To}.

\begin{theo}
Let $(N,g_0,J_0)$ be a compact K\"ahler-Einstein manifold with
nonpositive scalar curvature. Suppose all infinitesimal complex
deformations of $J_0$ are integrable. Then any deformation of
$g_0$ with constant scalar curvature must be K\"ahler-Einstein.
\end{theo}

This generalizes a result of \cite{ko} about Einstein
deformations.

\begin{theo}
Let $(N,g_0,J_0)$ be a compact K\"ahler-Einstein manifold with
negative scalar curvature. Suppose all infinitesimal complex
deformations of $J_0$ are integrable. Then there exists a
neighborhood $\mathcal{U}$ of $g_0$ in the space of smooth
Riemannian metrics on $N$ such that for any metric $g\in
\mathcal{U}$ with scalar curvature $S_g\geq S_{g_0}$
\begin{equation*}
{\rm Vol}(N,g)\geq {\rm Vol}(N, g_0)
\end{equation*}
and equality holds iff $g$ is a K\"ahler-Einstein metric with
negative scalar curvature.
\end{theo}

There are many examples satisfying the assumptions in the theorems
above. For example, the hypersurfaces of large enough degree in a
complex projective space. In fact we do not know any examples of
K\"ahler-Einstein manifolds of nonpositive scalar curvature which
do not satisfy the integrability condition for the complex
structure. It is likely that they all satisfy the integrability
condition, just as Calabi-Yau manifolds by virture of the
Bogomolov-Tian-Todorov theorem \cite{Bo}, \cite{t}, \cite{To}.

The study of the Yamabe constant, also called Schoen's $\sigma$
invariant, has attracted a lot of attention lately, Cf.
\cite{LB2}, \cite{BN}. This is motivated by a conjecture of Schoen
\cite{Sch1}, which says that the standard metric for manifolds
with constant sectional curvature realizes the Yamabe constant. In
other words, the standard metric is a global maximum for the
Yamabe invariant (they are called the supreme Einstein metrics in
\cite{LB2}). In view of the results of \cite{bcg2} for real
hyperbolic spaces and of \cite{LB1} for K\"ahler-Einstein
surfaces, it is tempting to conjecture the same for  more general
class of manifolds such as compact locallly symmetric spaces or
even K\"ahler-Einstein manifolds of negative scalar curvature.
Unfortunately, this is not true in higher dimensions, as a compact
simply connected manifold of dimension greater than or equal to
$5$ must have nonnegative Yamabe constant \cite{pe}, see also
\cite{S}.

There has been a lot of work recently concerning the stability of
Ricci flow \cite{gik}, \cite{sesum}, \cite{chen}, see also
\cite{chi}. The general question can be phrased as follows. If
$g_0$ is a metric such that the (renormalized) Ricci flow $g(t)$
starting from $g_0$ converges, is it true that the (renormalized)
Ricci flow $\tilde{g}(t)$ starting from all metrics $\tilde{g}_0$
that are sufficiently close to $g_0$ also converges? Using the
result of Natasa Sesum \cite{sesum}, we derive

\begin{theo} \label{ricciflow}
Let $(N,g_0,J_0)$ be a compact K\"ahler-Einstein manifold with
nonpositive scalar curvature. Suppose all infinitesimal complex
deformations of $J_0$ are integrable. Then the Ricci flow starting
from any Riemannian metric sufficient close to $g_0$ converges
exponentially to a K\"ahler-Einstein metric diffeomorphic to
$g_0$.
\end{theo}

The difference between this theorem and the well known result for
K\"ahler-Ricci flow on K\"ahler-Einstein manifolds with
nonpositive first Chern class \cite{cao} is that the Ricci flow
here starts with any metric nearby, rather than in a given K\"ahler
class. On the other hand, the result of \cite{cao} is a global
result in the sense that the initial metric is any metric in a given K\"ahler class.


This paper is organized as follows. We discuss spin$^c$ parallel
spinor and related Bochner type formula in the next section, and
prove the infinitesimal stability result. In Section 3, we discuss
the local stability results and applications.  In the final section,
relevant results from Kodaira-Spencer theory are recalled. We also
elaborate more on the examples and make some remarks.

{\em Acknowledgement:}  The authors wish to thank
Rick Schoen for stimulating discussions and encouragement.

\sect{Spin$^c$ parallel spinor and a Bochner type formula}

 We now assume $(M,g)$ is a compact Riemannian manifold
with a spin$^c$ structure. Thus, $w_2(M) \equiv c$, where $c\in
H^2(M, \mathbb Z)$ is the canonical class of the spin$^c$
structure. Let $\Sr^c\ra M$ denote the spin$^c$ spinor bundle and
$L \ra M$ the complex line bundle with $c_1(L)=c$. Then $\Sr^c=\Sr
\otimes L^{1/2}$, where the spinor bundle $\Sr$ may not exist
globally; similarly for the square root of the line bundle. An
excellent reference on spin geometry is Lawson and Michelsohn
\cite{LM}.

Let $E\ra M$ be a vector bundle with a connection. The curvature
is defined as
\begin{equation}
R_{XY}=-\grd_X\grd_Y+\grd_Y\grd_X+\grd_{[X,Y]}.
\end{equation}
If $M$ is a Riemannian manifold, then for the Levi-Civita
connection on $TM$, we have $R(X,Y,Z,W)=\lp R_{XY}Z,W\rp$. We
often work with an orthonormal frame $\{e_1,\ldots,e_n\}$ and its
dual frame $\{e^1,\ldots,e^n\}$. Set
$R_{ijkl}=R(e_i,e_j,e_k,e_l)$.

The spinor bundle $\Sr$, which may exist only locally, has a
natural connection induced by the Levi-Civita connection on $TM$.
For a spinor $\s$, we have
\begin{equation}
R_{XY}\s=\frac{1}{4}R(X,Y,e_i,e_j)e_ie_j\cdot \s.
\end{equation}
Given a unitary connection $\grd^L$ on $L$, we then obtain a
Clifford connection $\grd^c$ on $\Sr^c$. In fact,
$\grd^c=\grd\otimes 1 + 1\otimes \grd^{L^{1/2}}$ is the tensor
product connection for $\Sr^c$. Therefore, for a spin$^c$ spinor
$\s$, \be \label{curv}
R_{XY}\s=\frac{1}{4}R(X,Y,e_i,e_j)e_ie_j\cdot\s -\hf F(X, Y)
\s. \ee Here $F$ is the curvature form of $\grd^L$.

If $\s_0$ is a parallel spin$^c$ spinor, i.e., $\s_0$ is a
section of $\Sr^c$ such that $\grd^c_X \s_0=0$ for all $X$, then
$R_{XY}\s_0=0$. Hence we have
\begin{equation}\label{rxy=0}
R_{klij}e_ie_j\cdot \s_0 =2 F_{kl} \s_0.
\end{equation}

\begin{lem} If $\s_0$ is a parallel spin$^c$ spinor, then
\[ R_{kl}e_l\cdot \s_0 = F_{kl} e_l\cdot \s_0. \]
\end{lem}

\Pf From (\ref{rxy=0}) we have
\[ R_{klij}e_le_ie_j\cdot \s_0 =2 F_{kl} e_l\cdot \s_0. \]
But \ban  R_{klij} e_le_ie_j & = & \o{1}{3} \sum_{l,i,j \
distinct} (R_{klij} + R_{kijl} + R_{kjli}) e_le_ie_j  \\
& &  + \sum_{i,j} R_{kjij} e_je_ie_j + \sum_{i,j} R_{kiij}
e_ie_ie_j \\
& = & 2 R_{ki} e_i. \ean Here we have used the symmetries of
Riemann curvature tensor, including the first Bianchi identity.
Hence,
\[ R_{kl} e_l \cdot \s_0 =F_{kl} e_l \cdot \s_0 \]
as claimed. \qed

In the case that the spin$^c$ structure comes from a spin
structure, the line bundle $L$ is trivial; consequently $F=0$.
Thus $\Ric \equiv 0$ for manifolds with nonzero parallel spinor.

From now on, we assume $M$ has a parallel spin$^c$ spinor
$\s_0\not=0$, which, without loss of generality, is normalized to
be of unit length. We define, as in \cite{dww}, a linear map
$\Phi: S^2(M)\ra \Sr^c\ts T^*M$ by
\begin{equation} \label{susy}
\Phi(h)=h_{ij}e_i\cdot \s_0\ts e^j.
\end{equation}
It is easy to check that the definition is independent of the
choice of the orthonormal frame $\{e_1,\ldots,e_n\}$. The same
proof as in \cite{dww} again yields
\begin{lem} \label{phi}
The map $\Phi$ satisfies the following properties:
\begin{enumerate}
\item ${\rm Re}\,\lp \Phi(h),\Phi(\tilde{h})\rp=\lp h,\tilde{h}\rp
$, \item $\grd_X\Phi(h)=\Phi(\grd_Xh)$.
\end{enumerate} Here {\rm Re} denotes the real part.
\end{lem}

The following interesting Bochner type formula plays an important
role here.
\begin{lem}  \label{Bochner}
Let $h$ be a symmetric $2$-tensor on $M$. Then
\begin{equation} \label{bochner}
\Dirac^*\Dirac\Phi(h)=\Phi(\grd^*\grd h-2\cuv h - h\circ F + \Ric
\circ h).
\end{equation}
Here $(h\circ F)_{ij}=h_{ip}F_{pj}=-h_{ip}F_{jp}$ and $(\Ric\circ
h)_{ij}=R_{ip}h_{jp}$.
\end{lem}

\Rk Note that here we have implicitly extended our map $\Phi$ to
general (nonsymmetric) $2$-tensors with complex coefficients.

\Pf Choose an orthonormal frame $\{e_1,\ldots,e_n\}$ near a point
$p$ such that $\grd e_i=0$ at $p$. We compute at $p$, using Lemma
\ref{phi} and the Ricci identity,
\begin{align*}
\Dirac^*\Dirac\Phi(h) &=\grd_{e_k}\grd_{e_l}h(e_i,e_j)e_ke_le_i\cdot\s_0\ts e^j \\
&=-\grd_{e_k}\grd_{e_k}h(e_i,e_j)e_i\cdot\s_0\ts e^j
  -\frac{1}{2}R_{e_ke_l}h(e_i,e_j)e_ke_le_i\cdot \s_0\ts e^j\\
&=\Phi(\grd^*\grd h)
   +\frac{1}{2}R_{kljp}h_{ip}e_ke_le_i\cdot\s_0\ts e^j
   +\frac{1}{2}R_{klip}h_{pj}e_ke_le_i\cdot\s_0\ts e^j. \\
\end{align*}
By using twice the Clifford relation $e_ie_j+e_je_i=-2\delta_{ij}$
we have
\begin{align*}
\frac{1}{2}R_{kljp}h_{ip}e_ke_le_i\cdot\s_0
&=\frac{1}{2}R_{kljp}h_{ip}e_ie_ke_l\cdot\s_0+R_{kljp}h_{kp}e_l\cdot\s_0-R_{kljp}h_{lp}e_k\cdot
\s_0 \\
&=F_{jp}h_{ip}e_i\cdot\s_0 -2(\cuv h)_{kj}e_k\cdot\s_0.
\end{align*}
Here the last equality uses (\ref{rxy=0}). On the other hand,
\begin{equation*}
\frac{1}{2}R_{klip}h_{pj}e_ke_le_i\cdot\s_0=\hf h_{pj}
(R_{klip}e_ke_le_i\cdot\s_0)=R_{lp}h_{pj} e_l \cdot \s_0.
\end{equation*}
Putting these equations together we obtain our lemma. \qed

Once again, when the spin$^c$ structure comes from a spin
structure, the formula above becomes
\[ \Dirac^*\Dirac\Phi(h)=\Phi(\grd^*\grd h-2\cuv h),  \]
which recovers a formula of \cite{Wa2}, see also \cite{dww}.
By Lemma \ref{phi}, the stability result follows in this case \cite{dww}.

The existence of a parallel spin$^c$ spinor on a compact simply
connected manifold implies that the manifold is the product of a
K\"ahler manifold with a manifold with parallel spinor \cite{mo}.
Moreover, the spin$^c$ structure is the product of the canonical
spin$^c$ structure on the K\"ahler manifold with the spin
structure on the other factor.

We now assume that $(M, g)$ is a compact
K\"ahler manifold of real dimension $n=2m$.
Let $J$ be the parallel almost complex structure and
$\omega=g(J\cdot, \cdot)$ the K\"ahler form. The
complexified tangent bundle decomposes as
\begin{equation*}
TM  \otimes \mathbb C= T^{1,0}M \oplus T^{0,1}M.
\end{equation*}
The canonical spin$^c$ structure is given by
the anti-canonical line bundle
$L=K^{-1}=\Lambda^m\left(T^{1,0}(M)\right)$.
It has a canonical holomorphic connection induced from the
Levi-Civita Connection and the
curvature form $F=-\sqrt{-1}\rho$, where $\rho=\hbox{Ric}(J\cdot,\cdot)$
is the Ricci form.

The spinor bundle $\Sr^c(M)=\Sr^c_+ (M)\bigoplus\Sr^c_- (M)$ with
\begin{eqnarray*}
\Sr^c_{+}(M)&=&\bigoplus_{k\ even}\Lambda^{0,k}(M), \\
\Sr^c_{-}(M)&=&\bigoplus_{k\ odd} \Lambda^{0,k}(M).
\end{eqnarray*}
The Clifford multiplication is defined by
\[ v\cdot= \sqrt{2}(v^{0, 1}\!\wedge -v^{0,1}\lrcorner). \]
Here $v^{0,1}\lrcorner$ denotes the
contraction using the Hermitian metric.
The parallel spinor $\s\in C^\infty(\Sr^c_+(M))$ can be taken as the function
which is identically 1.

\medskip
\Rk
We would like to remark that the spin$^c$ structure and
related Bochner type formula are very useful in other context,
such as symplectic manifolds.
Given a symplectic manifold $(M,\omega)$ of dimension 2m,
we take an almost complex structure $J$ compatible with the
symplectic form $\omega$. This gives rise to a Riemannian metric $g=\omega(\cdot,J\cdot)$.
Then the formulation above in the K\"ahler case works perfectly well in
this generalized setting and defines a natural spin$^c$ structure on $M$.
The Levi-Civita connection $\nabla$ of $g$ induces a natural Hermitian connection $A$
on $\Lambda^m\left(T^{1,0}(M)\right)$ and hence a connection on the spinor bundle
$\Sr^c(M)$.
In general, the spinor $\s$ is not necessarily parallel.
In fact $\s$ is parallel iff $(M,\omega, J)$ is K\"ahler. However $\s$ is still a
harmonic spinor, a fact with several
interesting applications. Here we outline a simple example.
By the Lichenerowicz-Bochner formula we have
\[
\grd^*\grd
\s+\frac{S}{4}\s+\frac{1}{2}F_A\cdot\s=\Dirac^*\Dirac\s=0.
\]
Integrating by parts gives
\[
\int_M|\grd \s|^2+\frac{S}{4}+\frac{1}{2}\lp F_A\cdot\s,\s\rp =0.
\]
By straightforward calculations one can show
\[
|\grd \s|^2=\frac{1}{16}|\grd J|^2,
\]
\[
\lp F_A\cdot\s,\s\rp=-\frac{2\pi}{(m-1)!} C_1\wedge \omega^{m-1},
\]
where $C_1$ is the first Chern form.
Therefore we get the following interesting formula due to Blair \cite{b}
\[ \int_M\left(\frac{1}{4}|\grd J|^2+S\right)\frac{\omega^m}{m!}
=4\pi\int_M  C_1\wedge\frac{\omega^{m-1}}{(m-1)!}
\]
For more substantial applications of the spin$^c$ structure in symplectic geometry
we refer to the work of Taubes \cite{ta1, ta2}.

\bigskip
We now choose our orthonormal basis $e_1,
\cdots, e_{2m}$ so that $e_{m+i}=Je_i$. By a slight abuse of
notation, we denote $e_{\bi}=e_{m+i}=Je_i$. And similarly the
index $\bi$ denote $m+i$.
Hence, with $\s=1$ being the parallel
spin$^c$ spinor, we have
\begin{align} \lp \s, e_i e_j \cdot \s \rp = -\delta_{ij}, \ \ \ &
\lp \s, e_{\bi} e_{\bar{j}} \cdot \s \rp = -\delta_{ij}, \\
\lp \s, e_{\bi} e_j \cdot \s \rp = -\sqrt{-1}\delta_{ij}, \ \ \ &
\lp \s, e_{i} e_{\bar{j}} \cdot \s \rp =  \sqrt{-1}\delta_{ij}.
\end{align}

Now we compute
\begin{align*}
- \lp \Phi(h\circ F), \Phi(h) \rp =& \sum_{i, j, k, l, p=1}^{2m} F_{jp}
h_{ip} h_{kl} \lp e_i\cdot\s\ts e^j, e_k \cdot \s\ts e^l \rp \\
= -& \sum_{i, j, k, p=1}^{2m} F_{jp}
h_{ip} h_{kj} \lp \s, e_i e_k \cdot \s \rp \\
= & \sum_{j, p=1}^{2m} \sum_{i=1}^m  F_{jp} h_{ip} h_{ij} -
\sum_{j, p=1}^{2m} \sum_{i=1}^m  F_{jp} h_{ip} h_{\bi j}( \sqrt{-1}) \\
& - \sum_{j, p=1}^{2m} \sum_{i=1}^m F_{jp} h_{\bi p} h_{i j}(-\sqrt{-1})
+\sum_{j, p=1}^{2m} \sum_{i=1}^m F_{jp} h_{\bi p} h_{\bi j}.
\end{align*}

As the curvature of a unitary connection on a line bundle, $F$
is purely imaginary. Hence taking the real part (and using the
skew symmetry) yields: \be \label{rpe} - {\rm Re}\lp \Phi(h\circ F),
\Phi(h) \rp = -2\sqrt{-1}\sum_{j, p=1}^{2m} \sum_{i=1}^m  F_{jp}
h_{ip} h_{\bi j}. \ee

Similarly,
\be \label{rpe2}
{\rm Re}\lp \Phi(\Ric \circ h), \Phi(h) \rp = \sum_{i, j, p=1}^{2m}
  R_{ip} h_{pj} h_{i j}.
\ee

We are now ready to prove

\begin{theo} \label{inf}
If $(M, g_0)$ is a compact K\"ahler manifold with nonpositive
Ricci curvature, then $\grd^*\grd h-2\cuv h$ is positive
semi-definite on $S^2(M)$. That is,
\[ \lp \grd^*\grd h-2\cuv h, h \rp \geq \lp \Dirac \Phi(h), \Dirac \Phi(h) \rp \geq 0, \]
for any $h \in S^2(M)$. Moreover, in the case of negative Ricci curvature, $\grd^*\grd h-2\cuv h =0$ iff $\Dirac \Phi(h)= 0$
and $h$ is skew-hermitian.
\end{theo}

\Pf Since $\lp \Dirac \Phi(h), \Dirac \Phi(h) \rp \geq 0$, we have, by Lemmas~\ref{Bochner} and \ref{phi},
\begin{align*} 0 \leq & \lp \Dirac^*\Dirac \Phi(h), \Phi(h) \rp \\
= & {\rm Re} \lp \Phi(\grd^*\grd h-2\cuv h -h\circ F + \Ric
\circ h), \Phi(h) \rp \\
= & \lp \grd^*\grd h-2\cuv h, h \rp - {\rm Re} \lp \Phi(h\circ F), \Phi(h) \rp + {\rm Re} \lp
\Phi(\Ric \circ h), \Phi(h) \rp .
\end{align*}
That is,
\[ \lp \grd^*\grd h-2\cuv h, h \rp = \lp \Dirac \Phi(h), \Dirac \Phi(h) \rp - [ {\rm Re} \lp
\Phi(\Ric \circ h), \Phi(h) \rp - {\rm Re} \lp \Phi(h \circ F),
\Phi(h) \rp]. \]

For $L=K^{-1}$, as we remarked earlier, the curvature form $F=-\sqrt{-1}\rho$
where $\rho$ is the Ricci form. Since $g$ is K\"ahler, there is an orthonormal basis $e_1,
\cdots, e_{2m}$ such that $e_{m+i}=Je_i  (1\le i \le m)$ and the Ricci curvature is diagonal in this basis, i.e. $R_{ij} c_i \delta_{ij}  (1 \le i,j \le 2m)$ with $c_i = c_{m+i}$. Now
 \[
\rho (e_i,e_j) = {\rm Ric} (Je_i,e_j) = \left\{ \begin{array}{ll} c_i\delta_{m+i,j}, &  1\le i \le m \\
-c_i \delta_{m-i,j}, & m+1 \le i \le 2m \end{array} \right. .\]

It follows then from (\ref{rpe}) and (\ref{rpe2}) that
\[ -{\rm Re} \lp \Phi(h\circ F), \Phi(h) \rp = -2 \sum_{j=1}^m \sum_{i=1}^m c_j (h_{i\bar{j}}h_{\bi j}-h_{ij}h_{\bi \bar{j}}). \]
\[ {\rm Re} \lp \Phi(\Ric \circ h), \Phi(h) \rp = \sum_{i, j = 1}^{2m} c_i h^2_{ij}. \]
Hence,
\[ \lp \grd^*\grd h-2\cuv h, h \rp \geq  - [ \sum_{i, j = 1}^{2m} c_ih_{ij}^2
-2\sum_{j=1}^m \sum_{i=1}^m c_i (h_{i\bar{j}}h_{\bi j}-h_{ij}h_{\bi \bar{j}}) ]. \]
When $c_i \leq 0$ the right hand side is nonnegative
 by the Cauchy-Schwarz inequality.

If $c_i <0$, then $\grd^*\grd h-2\cuv h=0$ if and only if
$\Dirac\Phi(h)=0$ and
\[ h_{ij}=-h_{\bi \bar{j}},\quad h_{i\bar{j}}=h_{\bi j}. \]
That is, $h(JX, JY)=-h(X, Y)$. It follows that $h$ is  skew Hermitian.
\qed

Similarly, we have
\begin{theo} If $(M, g_0)$ is a compact K\"ahler manifold with nonnegative Ricci curvature, then the Lichnerowicz Laplacian
is positive
semi-definite on $S^2(M)$. That is,
\[ \lp \mathcal L_g h, h \rp \geq \lp \Dirac \Phi(h), \Dirac \Phi(h) \rp \geq 0, \]
for any $h \in S^2(M)$. Moreover, in the case positive Ricci curvature, $ \mathcal L_g h=0$
if and only if $\Dirac \Phi(h) = 0$ and $h$ is  Hermitian.
\end{theo}

\Pf The Lichnerowicz Laplacian is
\[ \mathcal L_gh=\grd^*\grd h-2\cuv h + \Ric \circ h + h\circ \Ric. \]
The same proof as above now goes through for $\mathcal L_g$. \qed

Note that the above computation in the case of K\"ahler-Einstein manifold with Einstein constant $c$
yields the following interesting Bochner-Lichnerowics-Weitzenbock formula:
\[   \lp \Dirac \Phi(h), \Dirac \Phi(h) \rp = \lp \grd^*\grd h-2\cuv h, h \rp + 2c\lp h_H, h_H \rp , \]
where $h_H$ denotes the Hermitian part of $h$. This unifies the two Weitzenbock formulas in
\cite[p. 362]{Besse}.

\sect{Local stability of K\"aher-Einstein metrics}

Theorem \ref{inf} says that for a K\"ahler-Einstein $(N,g_0)$ with
non-positive scalar curvature the operator \be\grd^*\grd h-2\cuv
h\ee is semi-positive definite on symmetric 2-tensors. A natural
and important question is to identify the kernal space
\be\label{ker} W_{g_0}=\{\ h\ |\ \op{tr}_{g_0}h=0,\d h=0,
\grd^*\grd h-2\cuv h=0\} \ee on the space of transverse traceless
symmetric 2-tensors. This is just the infinitesimal Einstein deformation space
studied in \cite{ko}. The case $c=0$ is essentially a Calabi-Yau
manifold, which has been studied with other manifolds admitting
parallel spinors in our previous paper \cite{dww}. We now focus on
the case $c<0$ using our approach.

By the proof of Theorem \ref{inf} \be W_{g_0}=\{\ h\ |\
\op{tr}_{g_0}h=0,\d h=0, h(J,J)=-h, \Dirac\Phi(h)=0\} \ee As
before we choose our orthonormal basis $e_1, \cdots, e_{2m}$ so
that $e_{m+i}=Je_i$.

Now for $1\leq i\leq m$ set
\[ X_i = \frac{e_i - \sqrt{-1} Je_i}{\sqrt{2}}, \ \ \
\bar{X}_i= \frac{e_i + \sqrt{-1}Je_i}{\sqrt{2}}. \] Then
$\{X_1,\ldots, X_m\}$ is a local unitary frame for $T^{1,0}M$ and
let its dual frame be $\{\theta^1,\ldots,\theta^m\}$. As $h$ is
skew-Hermition we have \be\label{skew}
h({X}_{i},\bar{X}_j)=h(\bar{X}_{j},{X}_i)=0 \ee By straightforward
computation we have for $h\in W=W_{g_0}$, \be
\Phi(h)=h(\bar{X}_{i},\bar{X}_j)\bar{\theta}^i\otimes\bar{\theta}^j.
\ee

This can be identified with
\be
\Psi(h)=h(\bar{X}_{i},\bar{X}_j)\bar{\theta}^i\otimes X_j\in \wedge^{0,1}(\Theta),
\ee
where $\Theta$ is the holomorphic tangent bundle.
We compute
\begin{eqnarray*}
\Dirac \Phi(h)&=&\sum_{k=1}^m
\left(\nabla_{e_k}h(\bar{X}_{i},\bar{X}_j)e_k\cdot\bar{\theta}^i\otimes\bar{\theta}^j
+\nabla_{e_{\bar{k}}}h(\bar{X}_{i},\bar{X}_j)e_{\bar{k}}\cdot\bar{\theta}^i\otimes\bar{\theta}^j\right) \\
&=&\sum_{k=1}^m
\left(\nabla_{e_k}h(\bar{X}_{i},\bar{X}_j)(\bar{\theta}^k\wedge\bar{\theta}^i-\d_{ik})\otimes\bar{\theta}^j
+\nabla_{e_{\bar{k}}}h(\bar{X}_{i},\bar{X}_j)
(\sqrt{-1}\bar{\theta}^k\wedge\bar{\theta}^i-\sqrt{-1}d_{ik})\otimes\bar{\theta}^j\right) \\
&=&\sqrt{2}\sum_{k=1}^m
\left(\nabla_{\bar{X}_k}h(\bar{X}_{i},\bar{X}_j)\bar{\theta}^k\wedge\bar{\theta}^i\otimes\bar{\theta}^j
-\nabla_{X_k}h(\bar{X}_{k},\bar{X}_j)\bar{\theta}^j\right)
\end{eqnarray*}
With $\Phi(h)$ identified as $\Psi(h)\in \wedge^{0,1}(\Theta)$,
the above calculation shows that the Dirac operator is then
$\sqrt{2}(\overline{\partial}-\overline{\partial}^*)$ on
$\Psi(h)$. Therefore $\Dirac\Phi(h)=0$ iff $\Psi(h)$ is harmonic.
On the other hand
\begin{eqnarray*}
\d h&=&\sum_{k=1}^{m}\left(\nabla_{e_k}h(e_k,\cdot)+\nabla_{e_{\bar{k}}}h(e_{\bar{k}},\cdot)\right) \\
&=&\sum_{k=1}^m (\nabla_{X_k}h(\bar{X}_k,\cdot)+\nabla_{\bar{X}_k} h (X_{k},\cdot))\\
&=&\sum_{j,k=1}^m\left(\nabla_{X_k}h(\bar{X}_{k},\bar{X}_j)\bar{\theta}^j
+\nabla_{\bar{X}_k}h(X_{k},{X}_j){\theta}^j\right)
\end{eqnarray*}
where in the last step we used (\ref{skew}).
This shows that $\d h=0$ automatically holds if $\Psi(h)$ is harmonic.
Therefore we have an injective homomorphism
\[
\Psi: W_{g_0}\ra H^1(N,\Theta).
\]
The image obviously consists of symmetric infinitesmial complex deformations.
To show that $\Psi$ is in fact onto we need to show all infinitesimal complex deformations are symmetric.
For this purpose we need a digression.

Let $N$ be a K\"ahler manifold with K\"ahler metric $\omega=\sqrt{-1}%
g_{i\overline{j}}dz^{i}\wedge d\overline{z}^{j}$. \ Given $\Psi=a_{\overline
{j}}^{i}d\overline{z}^{j}\otimes\frac{\partial}{\partial z^{i}}\in\wedge
^{0,1}(T^{1,0}N)$, we can consider the $(0,2)$-form%

\[
\psi=g_{i\overline{l}}a_{\overline{j}}^{i}d\overline{z}^{j}\wedge
d\overline{z}^{l}.
\]

\Rk One can work with a local unitary frame $\{X_1,\ldots, X_m\}$ and its dual frames
equally well, but it seems the calculations are easier working with local coordinates.

\medskip

We calculate%
\begin{align*}
\overline{\partial}\Psi &  =\frac{\partial a_{\overline{l}}^{i}}%
{\partial\overline{z}^{j}}d\overline{z}^{j}\wedge d\overline{z}^{l}%
\otimes\frac{\partial}{\partial z^{i}}\\
\overline{\partial}^{\ast}\Psi &  =-g^{k\overline{l}}\frac{\partial}%
{\partial\overline{z}^{l}}\rfloor\nabla_{\frac{\partial}{\partial z^{k}}}%
\Psi\\
&  =-g^{k\overline{l}}\left(  \frac{\partial a_{\overline{l}}^{i}}{\partial
z^{k}}+\Gamma_{kp}^{i}a_{\overline{l}}^{p}\right)  \frac{\partial}{\partial
z^{i}}%
\end{align*}
Suppose now that $\Psi$ is harmonic, i.e.
$\overline{\partial}\Psi=0, \overline{\partial}^{\ast}\Psi=0$. Then
we have
\begin{equation}\label{db=0}
\frac{\partial a_{\overline{l}}^{i}}
{\partial\overline{z}^{j}}=\frac{\partial a_{\overline{j}}^{i}}%
{\partial\overline{z}^{l}}
\end{equation}
and
\begin{equation}\label{db*=0}
g^{k\overline{l}}\left(  \frac{\partial a_{\overline{l}}^{i}}{\partial
z^{k}}+\Gamma_{kp}^{i}a_{\overline{l}}^{p}\right) =0
\end{equation}
Thus
\begin{equation*}
\overline{\partial}\psi =\left(g_{i\overline{l}}\frac{\partial a_{\overline{j}}^{i}}{\partial\overline{z}^{q}}
+a_{\overline{j}}^{i}\frac{\partial g_{i\overline{l}}}{\partial\overline{z}^{q}}\right)
d\overline{z}^{q}\wedge d\overline{z}^{j}\wedge d\overline{z}^{l} =0,
\end{equation*}
here we used \ref{db=0} and the fact that $\frac{\partial g_{i\overline{l}}}{\partial\overline{z}^{q}}$
is symmetric in $l$ and $q$.

We calculate
\begin{align*}
\overline{\partial}^{\ast}\psi &  =-g^{k\overline{l}}\frac{\partial}%
{\partial\overline{z}^{l}}\rfloor\nabla_{\frac{\partial}{\partial z^{k}}}%
\psi\\
&  =-g^{k\overline{l}}\frac{\partial}{\partial\overline{z}^{l}}\rfloor\left[
\left(  g_{i\overline{q}}\frac{\partial a_{\overline{j}}^{i}}{\partial z^{k}%
}+\frac{\partial g_{i\overline{q}}}{\partial z^{k}}a_{\overline{j}}%
^{i}\right)  d\overline{z}^{j}\wedge d\overline{z}^{q}\right]  \\
&  =-g^{k\overline{l}}\left[  \left(  g_{i\overline{q}}\frac{\partial
a_{\overline{l}}^{i}}{\partial z^{k}}+\frac{\partial g_{i\overline{q}}%
}{\partial z^{k}}a_{\overline{l}}^{i}\right)  d\overline{z}^{q}-\left(
g_{i\overline{l}}\frac{\partial a_{\overline{j}}^{i}}{\partial z^{k}}%
+\frac{\partial g_{k\overline{l}}}{\partial z^{i}}a_{\overline{j}}^{i}\right)
d\overline{z}^{j}\right]  \\
&  =-g_{i\overline{q}}g^{k\overline{l}}\left(  \frac{\partial a_{\overline{l}%
}^{i}}{\partial z^{k}}+\Gamma_{kp}^{i}a_{\overline{l}}^{p}\right)
d\overline{z}^{q}+\left(  \frac{\partial a_{\overline{j}}^{i}}{\partial z^{i}%
}+\frac{\partial\log\det G}{\partial z^{i}}a_{\overline{j}}^{i}\right)
d\overline{z}^{j}\\
&  =\left(  \frac{\partial a_{\overline{j}}^{i}}{\partial z^{i}}%
+\frac{\partial\log\det G}{\partial z^{i}}a_{\overline{j}}^{i}\right)
d\overline{z}^{j}%
\end{align*}
where in the last step we used \ref{db*=0}.
Therefore
\begin{align*}
\overline{\partial}\overline{\partial}^{\ast}\psi &  \left(  \frac{\partial^2 a_{\overline{j}}^{i}}{\partial z^{i}\partial\overline{z}^{l}}%
+\frac{\partial\log\det G}{\partial z^{i}}\frac{\partial a_{\overline{j}}^{i}}{\partial\overline{z}^{l}}
+\frac{\partial^{2}%
\log\det G}{\partial z^{i}\partial\overline{z}^{l}}a_{\overline{j}}^{i}\right)d\overline{z}^{l}\wedge d\overline{z}^{j}\\
&  =-\frac{\partial^{2}%
\log\det G}{\partial z^{i}\partial\overline{z}^{l}}a_{\overline{j}}^{i}d\overline{z}^{j}\wedge d\overline{z}^{l}
\end{align*}
where in the last step we used the fact that the first two coefficients are symmetric in $j$ and $l$ by (\ref{db=0}).
To summarize we have
\begin{equation}\label{fpsi}
\overline{\partial}\overline{\partial}^{\ast}\psi=R_{i\overline{l}}a_{\overline{j}}^{i}d\overline{z}^{j}\wedge d\overline
{z}^{l}%
\end{equation}

\bigskip

\begin{theo}\label{sym}
Let $(N,\omega_0)$ be a compact K\"ahler-Einstein manifold with negative scalar
curvature. Suppose $\Psi=a_{\overline{j}}^{i}d\overline{z}^{j}\otimes
\frac{\partial}{\partial z^{i}}\in\wedge^{0,1}(T^{1,0}N)$ is harmonic. Then
$\psi=g_{i\overline{l}}a_{\overline{j}}^{i}d\overline{z}^{j}\wedge
d\overline{z}^{l}=0$, i.e. $g_{i\overline{l}}a_{\overline{j}}^{i}$ is
symmetric in $i$ and $j$.
\end{theo}

\Pf
By the assumption we have $R_{i\overline{j}}=cg_{i\overline{j}}$ with $c<0$.
By (\ref{fpsi})
\[
\overline{\partial}\overline{\partial}^{\ast}\psi=c\psi.
\]
Therefore
\[
c\int_N |\psi|^2=\int_N\lp\overline{\partial}\overline{\partial}^{\ast}\psi,\psi\rp
=\int_N|\overline{\partial}^{\ast}\psi|^2
\]
Since $c<0$ we must have $\psi=0$.
\qed

\bigskip

\begin{Rk}
The discussion in Besse \cite{Besse}(12.96) contains some mistakes. The claim that skew-symmetric
infinitesimal complex deformations are in one-to-one correspondence with
holomorphic 2-forms is wrong. As the above calculation shows that in general
for a harmonc $\Psi\in H^{1}(N,\Theta)$ the corresponding $(0,2)$-form $\psi$
is not harmonic,
hence $\overline{\psi}$ is not holomorphic. The vanishing of the space of
skew-symmetric infinitesimal complex deformations on a K\"ahler-Einstein
manifold $N$ with negative
scalar curvature has nothing to do with the Hodge number $h^{2,0}=h^{0,2}$.
Take a compact complex hyperbolic surface $N$. By Calabi-Vesentini \cite{CV}
$H^{1}(N,\Theta)=0$. \ On the other hand, since the signature $\tau(N)>0$ and
the Euler number $\chi(N)=3$ $\tau(N)$, one can easily see by the Hodge index
theorem that $h^{2,0}(N)\neq0$ unless $N$ has the same Betti numbers as
$CP^{2}$, then a very special example constructed by Mumford. Therefore there
are compact complex hyperbolic surfaces $N$ with $H^{1}(N,\Theta)=0$ and
$h^{2,0}(N)\neq0$.
\end{Rk}

\medskip

\Rk In the Ricci flat case the above calculation shows that for a
harmonic $\Psi\in H^{1}(N,\Theta)$ the corresponding $(0,2)$-form
$\psi$ is indeed harmonic. Conversely one can show that a harmonic
$(0,2)$-form $\psi$ gives rise to a skew-symmetric infinitesimal
complex deformation. Therefore the space of skew-symmetric
infinitesimal complex deformations can be identified as the space
of holomorphic $(2,0)$-forms.

\medskip

As a corollary of Theorem \ref{sym} we now have a clear understanding of
the kernel $W_{g_0}$ of (\ref{ker}).
\begin{theo}
Let $(N,g_0,J_0)$ be a K\"ahler-Einstein manifold with negative scalar curvature. Then
$\Psi: W_{g_0}\ra H^1(N,\Theta)$ is an isomorphism.
\end{theo}

By the Kodaira-Spencer theorey, $H^1(N,\Theta)$ is the space of
infinitesimal complex deformations on $N$. In general these
infinitesimal deformations may not be integrable. But if they are
integrable, then the premoduli space of complex structures on $N$
is an manifold near $J$, with $H^1(N,\Theta)$ as the tangent space
(see the next section for more discussion). In this case we can
deduce various local results which we now explain. The argument is
by now standard, see \cite{dww} and Besson-Courtois-Gallot
\cite{bcg} where same type of results are established for Einstein
manifolds with negative sectional curvature.

We consider two well-known functionals, in addition to the first
eigenvalue $\lambda(g)$ considered in \cite{dww}. For a compact
manifold $(M,g)$ of dimension $n$ \be K(g)=\int_M|S_g|^{n/2}dV_g
\ee and \be Y(g)=\frac{\int_M S_gdV_g}{{\rm
Vol}(g)^{1-\frac{2}{n}}}. \ee

\begin{theo} \label{locforK}
Let $(N,g_0,J_0)$ be a compact K\"ahler-Einstein manifold with negative
scalar curvature. Suppose all infinitesimal complex deformations of $J_0$ are
integrable.
Then there
exists a neighborhood $\mathcal{U}$ of $g_0$ in the space of
smooth Riemannian metrics on $N$ such that
\begin{equation*}
\forall g\in \mathcal{U} \qquad K(g)\geq K(g_0)
\end{equation*}
and equality holds iff $g$ is a K\"ahler-Einstein metric with
negative scalar curvature. Moreover all Einstein metrics in
$\mathcal{U}$ are K\"ahler-Einstein with negative scalar
curvature.
\end{theo}

Since all infinitesimal complex deformations of $J_0$ are
integrable, the premoduli space of complex structures on $N$ is an
manifold near $J_0$, with $H^1(N,\Theta)$ as the tangent space. By
the uniqueness of K\"ahler-Einstein metric with negative scalar
curvature and the implicit function theorem, the moduli space
$\mathcal{E}$ of K\"ahler-Einstein metrics is an orbifold near
$g_0$, with $W_{g_0}\cong H^1(N,\Theta)$ as the tangent space.

Both functionals $K$ and $Y$ are scaling invariant, therefore we
can restrict ourselves to the space of Riemannian metric of volume
$1$, denoted by $\mathcal{M}$. By Ebin's slice theorem, there is a
real submanifold $\mathcal{S}$ containing $g_0$, which is a slice
for the action of the diffeomorphism group on $\mathcal{M}$. The
tangent space
\begin{equation}
T_{g_0}\mathcal{S}=\{h|\d_{g_0}h=0,
\int_N\op{tr}_{g_0}hdV_{g_0}=0.\}
\end{equation}
Let $\mathcal{C}\subset \mathcal{S}$ be the submanifold of constant scalar
curvatures metrics.

We need the following simple lemma from \cite{bcg}
\begin{lem}
Let $g$ be a metric with scalar curvature a negative constant and $g'$ a
metric conformal to $g$. Then $K(g')\geq K(g)$ and equality holds iff
$g'=g$.
\end{lem}

By this Lemma and the solution of the Yamabe problem, we only need
to prove $g_0$ is a local minimum for the functional $K$ on $\mathcal{C}$.
So it suffices to work on $\mathcal{C}$.
It is easy to see \be
T_{g_0}\mathcal{C}=\{h|\d_{g_0}h=0, \op{tr}_{g_0}h=0.\} \ee


Restricted on $\mathcal{C}$ and in a neighborhood of $g_0$, the
functional $K$ becomes $K(g)=|S_g|^m=(-S_g)^m$. Therefore to prove
that $g_0$ is a local minimum for $K$ on $\mathcal{C}$ is
equivalent to prove that $Y$ has a local maximum at $g_0$ on
$\mathcal{C}$. It is well known that $g_0$ is a critical point for
$Y$ and its Hessian at $g_0$ is given by \be
D^2Y(h,h)=-\frac{1}{2}\int_N\lp\grd^*\grd h-2\cuv h, h\rp. \ee

$\mathcal{C}$ contains the finite dimensional submanifold $\mathcal{E}$ of K\"ahler-Einstein
metrics, with tangent space $W_{g_0}\cong H^1(N,\Theta)$.
For any $g\in \mathcal{E}$, let $J$ be the associated complex
structure. We have
\[
\rho_g=\frac{S_g}{2m}\omega_g
\]
where $\omega_g$ is the associated K\"ahler form and $\rho_g$ the Ricci form.
Therefore
\[
Y(g)=-4m\pi\left(\frac{-C_1(N,J)^m[N]}{m!}\right)^{\frac{1}{m}},
\]
where $C_1(N,J)$ is the first Chern class of $(N,J)$. Thus $Y$ is
constant on $\mathcal{E}$.
Moreover, we have proved that $D^2Y$ is
negative definite on its normal bundle. Therefore there is a
possibly smaller neighborhood of $\mathcal{E}\subset \mathcal{C}$,
still denoted by $\mathcal{U}$, such that
\[
\forall g\in\mathcal{U}-\mathcal{E},\qquad Y(g)<Y(g_0).
\]
This proves Theorem \ref{locforK}.

We could have used $\lambda(g)$ instead of $Y(g)$ for the proof,
as in \cite{dww}.

\begin{theo}
Let $(N,g_0,J_0)$ be a compact K\"ahler-Einstein manifold with negative
scalar curvature. Suppose all infinitesimal complex deformations of $J_0$ are
integrable.
Then there
exists a neighborhood $\mathcal{U}$ of $g_0$ in the space of
smooth Riemannian metrics on $N$ such that for any metric
$g\in \mathcal{U}$ with scalar curvature $S_g\geq S_{g_0}$
\begin{equation*}
{\rm Vol}(N,g)\geq {\rm Vol}(N, g_0)
\end{equation*}
and equality holds iff $g$ is a K\"ahler-Einstein metric with
negative scalar curvature.
\end{theo}

\begin{Rk}
Though only a local result, it is quite remarkable to have volume
comparison under a lower bound for the scalar curvature. The
scalar curvature is a very weak geometric quantity and its effect
on a general Riemannian manifold $(M,g)$ of dimension $n$ can only
be detect infinitesimally by the following expansion for the
volume of a geodesic ball $B(p,r)$
\[
{\rm Vol}(B(p, r))=\omega_n r^n\left(1-\frac{S_g(p)}{6(n+2)}r^2+O(r^3)
\right) \quad\text{as } r\rightarrow 0
\]

\end{Rk}
\medskip

\Pf
We take the same $\mathcal{U}$ in Theorem \ref{locforK}. Then $\forall g\in
\mathcal{U}$ with $S_g\geq S_{g_0}$, we have $|S_g|^m\leq |S_{g_0}|^m$ since
$S_g<0$. Therefore $K(g)\leq |S_{g_0}|^m{\rm Vol}(N,g)$ while
$K(g_0)= |S_{g_0}|^m{\rm Vol}(N,g_0)$. The result then follows from
Theorem \ref{locforK}.
\qed

Theorem \ref{locforK} has another interesting interpretation.
Recall that the Yamabe invariant of a compact Riemannian manifold
$(M,g)$ of dimension $n$ is \be\label{yi} \mu(g) = \inf_{f\in\,
C^{\infty}(M),\ f>0} Y(f^{\frac{n-2}{4}}g)
\ee and it is a conformal invariant.
The Yamabe number of $M$ is defined as
\be
\sigma(M)={\rm sup}_{g} \mu (g).
\ee

We can now reformulate Theorem \ref{locforK}
as follows

\begin{theo} \label{locmax}
Let $(N,g_0,J_0)$ be a compact
K\"ahler-Einstein manifold with negative
scalar curvature. Suppose all infinitesimal complex deformations of $J_0$ are
integrable.
Then $g_0$ is a
local maximum of the Yamabe invariant.
\end{theo}

We end with a few remarks.
In \cite{Sch1} Schoen made the following conjecture
\begin{conj}
Let $(M,g_0)$ be a compact hyperbolic manifold. Then $\sigma(M)$ is
achieved by $g_0$ and only by $g_0$. In other words
\[ \forall g \quad \mu(g)\leq \mu(g_0) \]
and equality holds iff $g$ is conformal to $g_0$.
\end{conj}

It is also reasonable to make the same conjecture for other compact
locally symmetric spaces with negative sectional curvature.
In complex dimension 2 there have been some remarkable results proved
by LeBrun using Seiberg-Witten invariants. For example he proved that
on any compact K\"ahler-Einstein
surface with negative scalar curvature the Yamabe number is achieved
by the K\"ahler-Einstein metric. See \cite{LB1, LB2}.

In view of LeBrun's result and Theorem \ref{locmax} it is tempting
to extend the conjecture to K\"ahler-Einstein manifolds $(N,g_0)$
with negative scalar curvature in higher dimensions, namely that
$g_0$ should be a global maximum of the Yamabe invariant provided
all its infinitesimal complex deformations are integrable. But in
general this is not true. In fact any compact and simply connected
manifold $M$ of dimension $\geq 5$ has $\sigma(M)\geq 0$. This is
trivial if $M$ admits a metric of positive scalar curvature
\cite{S}. Otherwise it is proved by Petean \cite{pe}.

Another intriguing question is whether Theorem \ref{locforK} and
Theorem \ref{locmax} are still true if there are non-integrable
complex deformations. Then we have infinitesimal K\"aher-Einstein
deformations which can not be integrated to K\"aher-Einstein
metrics, but they may be integrated to Einstein metrics which are
not K\"ahler.

\sect{Kodaira-Spencer theory}

The deformation theory of complex structures was introduced by
Kodaira-Spencer in their seminal work \cite{ks1,ks2,ks3}. This
deep theory has played and still plays significant role in the
theory of complex manifolds. The relation between Kodaira-Spencer
theory and the deformation of Einstein metrics has been studied by
Koiso \cite{ko}. We review some relevant facts here in this
section and discuss some examples in more detail.

Let $M$ be a compact complex manifold and $\Theta$ the (sheaf of
germs of the) holomorphic tangent bundle of $M$. According to the
Kodaira-Spencer theory, the infinitesimal complex deformations are
described by the cohomology group $H^1(M, \Theta)$. For our
purpose, we are interested in the integrability of infinitesimal
complex deformations. Namely, when does every infinitesimal
deformation actually arise from a deformation of complex
structures? Let's recall first the so-called Theorem of Existence
in the Kodaira-Spencer theory \cite[Theorem 5.6]{Kod}.

\begin{theo}[Kodaira-Spencer] \label{kstoe} Let $M$ be a compact complex
manifold. If $H^2(M, \Theta)=0$, then there is a complex analytic
family with base $B$, $0\in B \subset \mathbb C^m$, such that the
fiber at $0$ is $M$ and the Kodaira-Spencer map at $0$ is an
isomorphism from $T_0B$ onto $H^1(M, \Theta)$.
\end{theo}

Recall that the Kodaira-Spencer map for a differentiable family of
compact complex manifolds assigns a tangent vector of the base to
the infinitesimal deformation along that direction. Thus the
condition $H^2(M, \Theta)=0$ implies that all infinitesimal
complex deformations are integrable.

We now discuss some examples from \cite{Kod}, where the reader is
referred to for complete detail.
\newline

\Ex 1). Blowups of $\mathbb{CP}^2$. For $M=\mathbb{CP}^2\#
k\overline{\mathbb{CP}^2}$ and $k\geq 5$, one has $H^0(M,
\Theta)=H^2(M, \Theta)=0$. Hence, all infintesimal deformations
are integrable. Incidentally, for $k \leq 4$, $H^1(M, \Theta)=0$.
Therefore the complex structure is rigid in these cases.

It is well-known that there exists K\"ahler-Einstein metrics on
$M$ if and only if $3 \leq k \leq 8$ by Tian's work \cite{tian}.
However, these K\"ahler-Einstein metrics have positive scalar
curvature. Hence our results do not apply. In fact, other than $\mathbb{CP}^2$
itself, these are unstable, \cite{chi}.
\newline

\Ex 2). Surfaces of arbitrary degree. For a non-singular surface
$M$ of degree $h$ in $\mathbb{CP}^3$, one has \[ \dim H^2(M,
\Theta)=\hf (h-2)(h-3)(h-5). \] Thus, $H^2(M, \Theta)=0$ for $h=2,
3, 5$. Since $c_1(M)=(4-h)H$ where $H$ is the hyperplane class,
$M$ has K\"ahler-Einstein metrics with negative scalar curvature
if $h\geq 5$, by the Calabi-Aubin-Yau Theorem. Hence our results
apply to the non-singular surface of degree $5$.

As one can see here, in general, the condition $H^2(M, \Theta)=0$,
which guarantees the integrability of all infinitesimal complex
deformations, is very restrictive. Indeed, there are many examples
which do not satisfy this condition but still, all their
infinitesimal complex deformations are integrable. In fact,
understanding the reason behind this is one of the motivations for
Kodaira-Spencer.

From the Kodaira-Spencer theory, if an infinitesimal complex
deformation $\theta \in H^1(M, \Theta)$ is integrable, then
\[ [\theta, \theta]=0. \]
This is in fact the first order obstruction. One thus expects that
there should be non-integrable infinitesimal deformations.
However, the cohomology group $H^1(M, \Theta)$ turns out to be
surprisingly difficult to compute. And in the many examples where
it can be computed, the infinitesimal deformations turn out to be
integrable.

Recall that a complex analytic family of compact complex manifolds
is said to be effective (or minimal) if its Kodaira-Spencer map is
injective. It is called a complete (or versal) family if every
other (sufficiently small) family of deformations can be induced
from this family via pullback of a holomorphic map. Now whenever
there is an effective complete family with base $B \ni 0$ a domain
in $\mathbb C^m$, such that the fiber at $0$ is $M$,
Kodaira-Spencer defines the number of moduli $m(M)=m$ to be the
dimension of the base. Then the question of whether all
infinitesimal deformations are integrable can be reinterpreted as
when the equality $m(M)=\dim H^1(M, \Theta)$ holds, of which
Kodaira-Spencer refers as the fundamental guiding question in the
Kodaira-Spencer theory.

By the Theorem of Completeness in the Kodaira-Spencer theory
\cite[Theorem 6.1]{Kod}, which says that a complex analytic family
of compact complex manifolds with surjective Kodaira-Spencer map
is complete, the complex analytic family in Theorem \ref{kstoe} is
an effective complete family with base dimension $\dim H^1(M,
\Theta)$.

More generally, if $M$ is a compact complex manifold for which
there is a complex analytic family of deformations whose
Kodaira-Spencer map is surjective, then all infinitesimal
deformations are integrable. This is the case where many examples
can be found. (We will see the converse in a moment.)
\newline

\Ex 3). Hypersurfaces in $\mathbb{CP}^n$. If $M$ is a hypersurface
in $\mathbb{CP}^n$ of degree $d$, one can construct a complex
analytic family of deformations of $M$ by varying the coefficients
of the defining equation of $M$. This family has surjective
Kodaira-Spencer map.  There
are many examples in this class which admits K\"ahler-Einstein
metrics with negative scalar curvature.  Let $N\subset {\bf
CP}^{m+1}$ be a smooth algebraic hypersurface of degree $d>m+2$.
Then the 1st Chern class $c_1(N)<0$. By the theorem of
Calabi-Aubin-Yau, there is a K\"ahler-Einstein metric with
negative scalar curvature $g_0$ on $N$. It is shown in
\cite[p219]{Kod} that
\[
{\rm dim}H^1(N,\Theta)={m+1+d \choose d}-(m+2)^2.
\]

Going back to the question of existence, without any assumptions,
there is the Kuranishi Theorem \cite{ku}.

\begin{theo}[Kuranishi] For any compact complex manifold $M$,
there exists a complete complex analytic family with base $B$,
$0\in B$ such that the fiber at $0$ is $M$. Moreover $B$ is a
complex analytic subset of $\mathbb C^m$, where $m=\dim H^1(M,
\Theta)$, defined by $l$ holomorphic equations, with $l=\dim
H^2(M, \Theta)$.
\end{theo}

It can be deduced from Kuranishi's theorem that if every
infinitesimal complex deformations are integrable, then the
Kuranishi family above is a complex analytic family whose
Kodaira-Spencer map is an isomorphism (and hence surjective).
Thus, our integrability assumption is equivalent to the existence
of complex analytic family of deformations whose Kodaira-Spencer
map is surjective.

As pointed out by Koiso \cite{ko}, the existence of such a family
has significant implication for the moduli space of
(K\"ahler-)Einstein metrics.

\begin{theo} Let $M$ be a compact complex manifold and $(g, J)$ be
a K\"ahler-Einstein structure on $M$. Assume that the complex
structure $J$ belongs to a complex analytic family of complex
structures with surjective Kodaira-Spencer map. Moreover, if the
scalar curvature is positive, assume further that there is no
nonzero hermitian infinitesimal Einstein deformations and also no
nonzero holomorphic vector field. Then the local premoduli space
of Einstein metrics around $g$ is a manifold with tangent space at
$g$ the space of infinitesimal Einstein deformations. Moreover,
any Einstein metric in it is K\"ahler (with respect to some
complex structure).
\end{theo}

It should be pointed out that there are indeed examples of compact
complex manifolds with non-integrable infinitesimal complex
deformations \cite[p319]{Kod}. However, we do not know any
examples of K\"ahler-Einstein manifolds with negative scalar
curvature which does not satisfy the integrability condition. In
view of the Bogmolov-Tian-Todorov theorem \cite{Bo,t,To} in the
Calabi-Yau case, we have the following very interesting question.

{\bf Question:} Is it true that on any compact K\"ahler-Einstein
manifolds with negative scalar curvature, the universal
deformation space of complex structures is smooth?


\begin{thebibliography}{SoWei2}


\bibitem[Bes87]{Besse}
Arthur~L. Besse.
\newblock {\em Einstein manifolds}.
\newblock Springer-Verlag, Berlin, 1987.

\bibitem[BCG91]{bcg}
G. Besson, G. Courtois, S. Gallot, {\it Volume et entropie minimale des espaces localement
sym\'etriques}, Invent. math. 103, 417-445 (1991).

\bibitem[BCG95]{bcg2}G. Besson, G. Courtois, S. Gallot, {\it Entropies et rigiditŽs des espaces localement symŽtriques de courbure strictement nŽgative}, Geom. Funct. Anal. 5 (1995), no. 5, 731--799.

\bibitem[Bl92]{b} D. Blair, {\em The ``total scalar curvature" as a symplectic invariant and related results }.
Proceedings of the 3rd Congress of Geometry (Thessaloniki, 1991), 79--83,
Aristotle Univ. Thessaloniki, Thessaloniki, 1992.

\bibitem[Bo78]{Bo} F. A. Bogomolov, {\em Hamiltonian K\"ahler manifolds}, Dolk. Akad. Nauk SSSR
243(1978), no. 5,  1101-1104.

\bibitem[B\"o05]{Boehm} C. B\"ohm, {\em Unstable Einstein metrics},
Math. Z. (2005).

\bibitem[BN04]{BN} H. Bray, A.  Neves, {\em Classification of prime 3-manifolds with Yamabe invariant greater than
$\mathbb{RP}\sp 3$},  Ann. of Math. (2)  159  (2004),  no. 1, 407--424.

\bibitem[CV60]{CV} E. Calabi, E. Vesentini, {\em On compact, locally symmetric K\"ahler manifolds},
Ann. of Math. (2) 71 (1960) 472-507.

\bibitem[Cao85]{cao} H. Cao, {\em Deformation of K\"ahler metrics to K\"ahler-Einstein metrics
 on compact K\"ahler manifolds}, Invent. Math. 81 (1985), no. 2, 359--372.

\bibitem[CHI04]{chi} Huai-Dong Cao, Richard S. Hamilton, Tom Ilmanen,
{\it Gaussian densities and stability for some Ricci solitons},
math.DG/0404165.

\bibitem[Ch05]{chen} X. Chen, {\em On the lower bound of energy functional $E_1$ (I)-- a stability
theorem on the Kaehler Ricci flow}, math.DG/0502196.


\bibitem[DWW04]{dww} X. Dai, G. Wei, X. Wang, {\it On the stability of Riemannian manifold
with parallel spinors}, Inventiones Mathematicae (2005).

\bibitem[GIK02]{gik} C. Guenther, J. Isenberg, D. Knopf, {\em
Stability of the Ricci flow at Ricci-flat metrics}, Comm. Anal.
Geom. 10 (2002), no. 4, 741--777.

\bibitem[H74]{H} N. Hitchin, {\em Harmonic spinors}, Adv. in Math.
14(1974), 1-55.



\bibitem[J00]{J} D. Joyce, {\em Compact manifolds with special holonomy},
Oxford Univ. Prss, Oxford, 2000.

\bibitem[Kod86]{Kod} K. Kodaira, {\em Complex manifolds and deformations of complex
structures}, Springer-Verlag, 1986.

\bibitem[KS58-1]{ks1} K. Kodaira and D. Spencer, {\em On deformations of complex
analytic structures}, I-II, Ann. Math., 67(1958), 328-466.

\bibitem[KS58-2]{ks2} K. Kodaira and D. Spencer, {\em A theorem of completeness for complex
analytic fibre spaces}, Acta Math., 100(1958), 281-294.

\bibitem[KS60]{ks3} K. Kodaira and D. Spencer, {\em On deformations of complex
analytic structures}, III, Ann. Math., 71(1960), 43-76.

\bibitem[Ko79]{ko1} N. Koiso, {\em A decomposition of the space ${\cal M}$ of Riemannian metrics on a manifold},
  Osaka J. Math.  16  (1979), no. 2, 423--429.

\bibitem[Ko80]{ko2} N. Koiso, {\em Rigidity and stability of Einstein metrics---the case of compact symmetric spaces},
  Osaka J. Math.  17  (1980), 51--73.

\bibitem[Ko83]{ko} N. Koiso, {\em Einstein metrics and complex structures}, Invent.
Math. 73(1983), 71-106.

\bibitem[Ku62]{ku} M. Kuranishi, {\em On the locally complete families of complex
analytic structures}, Ann. Math., 75(1962), 536-577.


\bibitem[LM89]{LM}
H.~Blaine Lawson, Jr. and Marie-Louise Michelsohn.
\newblock {\em Spin geometry}.
\newblock Princeton University Press, Princeton, NJ, 1989.

\bibitem[Le95]{LB1} LeBrun
{\it Einstein metrics and Mostow rigidity}, MRL 2(1995) 1-8.

\bibitem[Le99]{LB2} LeBrun {\it Einstein metrics and the Yamabe problem}, Trends in
mathematical physics, 353-376, AMS, Providence, RI, 1999.

\bibitem[L63]{L} A. Lichnerowicz, {\em Spineurs harmonique}, C.
R. Acad. Sci. Paris, S\'er. A-B,257(1963), 7-9.

\bibitem[M96]{mor} J. Morgan, {\em The Seiberg-Witten equations and applications
to the topology of smooth four-manifolds}, Princeton Univ. Press,
Princeton, NJ, 1996.

\bibitem[Mo97]{mo} A. Moroianu, {\em Parallel and Killing Spinors on Spinc Manifolds},
Commun. Math. Phys. 187, 417-427 (1997).

\bibitem[Pe00]{pe} J. Petean, {\em The Yamabe invariant of simply connected manifolds},
J. Reine Angew. Math. 523 (2000), 225--231.


\bibitem[P02]{P} G. Perelman. {\em The entropy formula for the Ricci flow and
its geometric
applications}, math.DG/0211159.


\bibitem[Sch84]{Sch2}
Richard~M. Schoen. {\em Conformal deformation of a Riemannian
metric to constant scalar curvature}, J. Differential Geom. 20
(1984), no. 2, 479--495.

\bibitem[Sch89]{Sch1}
Richard~M. Schoen.
\newblock Variational theory for the total scalar curvature functional for
  {R}iemannian metrics and related topics.
\newblock In {\em Topics in calculus of variations (Montecatini Terme, 1987)},
  pages 120--154. Springer, Berlin, 1989.



\bibitem[Se04]{sesum} N. Sesum, {\em Linear and dynamical stability of Ricci
flat metrics}, math.DG/0410062.

\bibitem[S92]{S} S. Stolz, {\em Simply connected manifolds of
positive scalar curvature}, Ann. Math., 136(1992), 511-540.

\bibitem[Ta94]{ta1} C. Taubes, {\em The Seiberg-Witten invariants and symplectic forms},
Math. Res. Lett. 1 (1994), no. 6, 809--822.
\bibitem[Ta95]{ta2} C. Taubes, {\em More constraints on symplectic forms from Seiberg-Witten invariants},
Math. Res. Lett. 2 (1995), no. 1, 9--13.


\bibitem[T86]{t} G. Tian
\newblock Smoothness of the universal deformation space of compact {Calabi-Yau} manifolds and its {Petersson-Weil} metric.
\newblock In {\em Mathematical aspects of string theory}, pages
629-646, World Scientific, 1986.

\bibitem[T97]{tian} G. Tian, {\em K\"ahler-Einstein metrics with positive scalar curvature},
Invent. Math. 130(1997), no. 1, 1--37.

\bibitem[To89]{To} A. Todorov, {\em The Weil-Petersson geometry of the
moduli space of ${\rm SU}(n\geq 3)$ (Calabi-Yau) manifolds. I.},
Comm. Math. Phys. 126 (1989), no. 2, 325--346.

\bibitem[Wa89]{Wa} M. Wang, {\it Parallel spinors and parallel forms},
Ann. Global Anal. Geom., 7(1989), no.1, 59-68.

\bibitem[Wa91]{Wa2} M. Wang, {\it Preserving parallel spinors under metric deformations},
Indiana Univ. Math. J., 40 (1991), no. 3, 815--844.

\bibitem[Y77]{y1} S. T. Yau, {\em Calabi's conjecture and some new results in algebraic
geometry},  Proc. Nat. Acad. Sci. U.S.A. 74 (1977), no. 5,
1798--1799.

\bibitem[Y78]{y2} S. T. Yau, {\em On the Ricci curvature of a compact K\"ahler manifold
and the complex Monge-Amp\`ere equation. I.}, Comm. Pure Appl.
Math. 31 (1978), no. 3, 339--411.

\bibitem[Ye93]{ye} R. Ye, {\em Ricci flow, Einstein metrics and space forms},
 Trans. Amer. Math. Soc. 338 (1993), no. 2, 871--896.

\end{thebibliography}
\end{document}